\documentclass[12pt,reqno]{amsart}

\usepackage{fullpage,amsmath,amssymb,amsfonts,mathrsfs,bm,graphicx,hyperref}

\def\bbone{{\mathchoice {\rm 1\mskip-4mu l} {\rm 1\mskip-4mu l}
{\rm 1\mskip-4.5mu l} {\rm 1\mskip-5mu l}}}

\newtheorem{theorem}{Theorem}[section]

\newtheorem{lemma}{Lemma}[section]
\newtheorem{proposition}{Proposition}[section]
\newtheorem{conjecture}{Conjecture}[section]

\begin{document}

\author{Abdelmalek Abdesselam}
\address{Abdelmalek Abdesselam, Department of Mathematics,
P. O. Box 400137,
University of Virginia,
Charlottesville, VA 22904-4137, USA}
\email{malek@virginia.edu}

\author{Pedro Brunialti}
\address{Pedro Brunialti, Department of Mathematics,
P. O. Box 400137,
University of Virginia,
Charlottesville, VA 22904-4137, USA}
\email{eyu4dh@virginia.edu}

\author{Tristan Doan}
\address{Tristan Doan, Department of Mathematics,
P. O. Box 400137,
University of Virginia,
Charlottesville, VA 22904-4137, USA}
\email{tgj9xx@virginia.edu}

\author{Philip Velie}
\address{Philip Velie, Department of Mathematics,
P. O. Box 400137,
University of Virginia,
Charlottesville, VA 22904-4137, USA}
\email{pmv8ev@virginia.edu}

\title{A bijection for tuples of commuting permutations and a log-concavity conjecture}

\begin{abstract}
Let $A(\ell,n,k)$ denote the number of $\ell$-tuples of commuting permutations of $n$ elements whose permutation action results in exactly $k$ orbits or connected components.
We provide a new proof of an explicit formula for $A(\ell,n,k)$ which is essentially due to Bryan and Fulman, in their work on orbifold higher equivariant Euler characteristics. Our proof is self-contained, elementary, and relies on the construction of an explicit bijection, in order to perform the $\ell+1\rightarrow \ell$ reduction. 
We also investigate a conjecture by the first author, regarding the log-concavity of $A(\ell,n,k)$ with respect to $k$. The conjecture generalizes a previous one by Heim and Neuhauser related to the Nekrasov-Okounkov formula. 
\end{abstract}

\maketitle


\section{Introduction}

For $n\ge 0$, let us denote by $[n]$ the finite set $\{1,\ldots,n\}$, and by $\mathfrak{S}_n$
the symmetric group of permutations of $[n]$. For $\ell\ge 0$, we consider the set of ordered $\ell$-tuples of commuting permutations
\[
\mathscr{C}_{\ell,n}:=\left\{\ 
(\sigma_1,\ldots,\sigma_{\ell})\in (\mathfrak{S}_n)^{\ell}\ |\ \forall i,j,\  \sigma_i\sigma_j=\sigma_j\sigma_i
\ \right\}\ .
\]
For a tuple $(\sigma_1,\ldots,\sigma_{\ell})$ of (non-necessarily commuting) permutations, let $\langle \sigma_1,\ldots,\sigma_{\ell}\rangle$ 
be the subgroup they generate inside $\mathfrak{S}_n$. The obvious action of $\mathfrak{S}_n$ on $[n]$ restricts to an action of $\langle \sigma_1,\ldots,\sigma_{\ell}\rangle$
with a number of orbits which we will denote by $\kappa( \sigma_1,\ldots,\sigma_{\ell})$.
For $0\le k\le n$, we let $\mathscr{C}_{\ell,n,k}$ be the subset of $\mathscr{C}_{\ell,n}$ made of tuples for which $\kappa( \sigma_1,\ldots,\sigma_{\ell})=k$. We finally define our main object of study
\[
A(\ell,n,k):=|\mathscr{C}_{\ell,n,k}|\ ,
\]
where, as usual, $|\cdot|$ denotes the cardinality of finite sets.
Our main result is a new proof of the following theorem giving an explicit, albeit complicated, formula for the $A(\ell,n,k)$.

\begin{theorem}
We have
\[
A(\ell,n,k)=\frac{n!}{k!}\times
\sum_{n_1,\ldots,n_k\ge 1}
\bbone\{n_1+\cdots+n_k=n\}\times \prod_{i=1}^{k}
\left[\frac{B(\ell,n_i)}{n_i}\right]\ ,
\]
where $\bbone\{\cdots\}$ denotes the indicator function of the condition between braces, and $B(\ell,\cdot)$ is the multiplicative function (in the number theory sense, i.e., $B(\ell,ab)=B(\ell,a)B(\ell,b)$ when $a,b$ are coprime) which satisfies
\[
B(\ell,q^m)=\frac{(q^{\ell}-1)(q^{\ell+1}-1)\cdots(q^{\ell+m-1}-1)}{(q-1)(q^2-1)\cdots(q^m-1)}\ ,
\]
when $m\ge 0$ and $q$ is a prime number.
\label{mainthm}
\end{theorem}

Our motivation for considering the previous theorem is the following log-concavity conjecture
by the first author.

\begin{conjecture} (Abdesselam~\cite{Abdesselam2})
For all $\ell\ge 1$, all $n\ge 3$, and for all $k$ such that $2\le k\le n-1$,
\[
A(\ell,n,k)^2\ge A(\ell,n,k-1)\ A(\ell,n,k+1)\ .
\]
\label{mainconj}
\end{conjecture}

The case $\ell=1$, included for esthetic coherence, is not conjectural. Since $A(1,n,k)=c(n,k)$, the unsigned Stirling number of the first kind, the stated log-concavity
property is a well known fact (see, e.g.,~\cite{Abdesselam1} and references therein).
The case $\ell=2$ is a conjecture by Heim and Neuhauser~\cite{HeimN} related to the Nekrasov-Okounkov formula~\cite{NekrasovO,Westbury}, as will be explained in \S\ref{conjsec}. The case ``$\ell=\infty$'' is proved in~\cite{Abdesselam2}.
The form in which Theorem \ref{mainthm} is stated is the one needed for the proof given in~\cite{Abdesselam2}, and we did not see this precise formulation in the literature.
However, we do not claim Theorem \ref{mainthm} is new. Indeed, it follows easily from the following identity by Bryan and Fulman~\cite{BryanF}
\begin{equation}
\sum_{n=0}^{\infty}\sum_{k=0}^{n}\frac{1}{n!}\ A(\ell,n,k)\ x^k u^n=
\prod_{d_1,\ldots,d_{\ell-1}=1}^{\infty}
(1-u^{d_1\cdots d_{\ell-1}})^{-x\, d_1^{\ell-2} d_2^{\ell-3}\cdots d_{\ell-2}}\ ,
\label{BFidentity}
\end{equation}
which holds in the ring of formal power series $\mathbb{C}[[x,u]]$.
To see how Theorem \ref{mainthm} can be derived from (\ref{BFidentity}), first (re)define, for $\ell\ge 1$ and $n\ge 1$,
\begin{equation}
B(\ell,n):=\sum_{s_1|s_2|\cdots|s_{\ell-1}|n}s_1\cdots s_{\ell-1}\ ,
\label{altdefeq}
\end{equation} 
where the sum is over tuples of integers $s_1,\ldots,s_{\ell-1}\ge 1$
which form an ``arithmetic flag'', namely, such that $s_1$ divides $s_2$,  $s_2$ divides $s_3$,\ldots, $s_{\ell-1}$ divides $n$.
In particular, $B(1,n)=1$, and $B(2,n)=\sigma(n)$ the divisor sum from number theory.
Since the divisor lattice factorizes over the primes, it is clear from the alternative definition (\ref{altdefeq}), that the $B(\ell,\cdot)$ is a mutiplicative function, in the number theory sense.
Its computation reduces to the prime power case. If $q$ is a prime and $m\ge 0$, then we have
\begin{eqnarray*}
B(\ell,q^m) &=& \sum_{0\le m_1\le\cdots\le m_{\ell-1}\le m}\ q^{m_1+\cdots+m_{\ell-1}} \\
 & = & \sum_{\lambda\subset (m)^{\ell-1}} q^{|\lambda|}\\
& = & \left[
\begin{array}{c} 
m+\ell-1 \\
m
\end{array}
\right]_q \\
 & = & \frac{(q^{\ell}-1)(q^{\ell+1}-1)\cdots(q^{\ell+m-1}-1)}{(q-1)(q^2-1)\cdots(q^m-1)}\ .
\end{eqnarray*}
Here, we changed variables to the integer partition $\lambda=(m_{\ell-1},m_{\ell-2},\ldots,m_{1})$ with weight $|\lambda|$ and whose shape is contained in the rectangular partition $(m)^{\ell-1}$ with $\ell-1$ parts equal to $m$. Finally, we used the well known formula for the sum over $\lambda$ as a Gaussian polynomial or $q$-binomial coefficient (see, e.g.,~\cite[Prop. 1.7.3]{Stanley}).
This shows the equivalence between (\ref{altdefeq}) and the definition given in Theorem \ref{mainthm}. By changing variables from $s_1,\ldots,s_{\ell-1}$ to $d_1,\ldots,d_{\ell}$ given by
\[
d_1=s_1\ ,\ d_2=\frac{s_2}{s_1}\ ,\ldots,\ d_{\ell-1}=\frac{s_{\ell-1}}{s_{\ell-2}}\ ,\  d_{\ell}=\frac{n}{s_{\ell-1}}\ ,
\]
we can also write
\[
B(\ell,n)=\sum_{d_1\cdots d_{\ell}=n} d_1^{\ell-1}d_{2}^{\ell-2}\cdots d_{\ell-1}\ ,
\]
as a multiple Dirichlet convolution of power functions (see, e.g.,~\cite{Moller} where the connection to $q$-binomial coefficients was also noted). The last formula is also consistent with the extreme $\ell=0$ case, where $B(0,n)=\bbone\{n=1\}$.

We then have the following easy formal power series computations
\begin{eqnarray*}
\sum_{n=1}^{\infty}\frac{B(\ell,n)}{n}\ u^n &= &
\sum_{d_1,\ldots,d_{\ell}\ge 1} \frac{d_1^{\ell-1}d_{2}^{\ell-2}\cdots d_{\ell-1}}{d_1\cdots d_\ell} \times u^{d_1\cdots d_{\ell}}\\
 & = & \sum_{m\ge 1} B(\ell-1,m)\times
\sum_{d_{\ell}\ge 1}\frac{(u^m)^{d_{\ell}}}{d_{\ell}}\\
 & = & \sum_{m\ge 1} B(\ell-1,m) \times\left(-\log(1-u^m)\right)\ ,
\end{eqnarray*}
where we introduced the new summation index $m:=d_1\cdots d_{\ell-1}$.
Multiplying by $x$, and taking exponentials gives
\begin{equation}
\exp\left(x \sum_{n=1}^{\infty}\frac{B(\ell,n)}{n}\ u^n\right)=
\prod_{m=1}^{\infty} (1-u^m)^{-xB(\ell-1,m)}\ ,
\end{equation}
which is the right-hand side of (\ref{BFidentity}) when collecting factors according to $m:=d_1\cdots d_{\ell-1}$.
We have thus shown that (\ref{BFidentity}) can be rewritten as
\begin{equation}
\sum_{n=0}^{\infty}\sum_{k=0}^{n}\frac{1}{n!}\ A(\ell,n,k)\ x^k u^n=
\exp\left(x \sum_{n=1}^{\infty}\frac{B(\ell,n)}{n}\ u^n\right)\ .
\end{equation}
Extracting coefficients of monomials in $x$ and $u$, on both sides, immediately yields Theorem \ref{mainthm}.
In the article~\cite{BryanF}, $x$ is assumed to be the Euler characteristic of a manifold. However, their proof of (\ref{BFidentity}) holds if $x$ merely is a formal variable.
Their work was aiming at generalizing the ``stringy'' orbifold Euler characteristic~\cite{DixonHVW,AtiyahS}, from sums over pairs of commuting permutations, to commuting tuples of arbitrary length $\ell$. Another motivation for their work was the study by Hopkins, Kuhn, and Ravenel~\cite{HopkinsKR} of a hierarchy of cohomology theories where the $\ell$-th level seemed to crucially involve $\ell$-tuples of commuting elements of a finite group such as $\mathfrak{S}_n$. 
The group-theoretic proof by Bryan and Fulman involved a delicate analysis of conjugacy classes in wreath products. Another proof one can find in the literature is the one by White~\cite{White}.
It uses the remark that $\mathscr{C}_{\ell,n}$ is in bijection with
${\rm Hom}(\mathbb{Z}^{\ell},\mathfrak{S}_n)$, namely, the set of group homomorphisms
from the additive group $\mathbb{Z}^\ell$ to the symmetric group $\mathfrak{S}_n$, i.e.,
$\mathbb{Z}^{\ell}$ actions on a set of $n$ elements.
The proof by White also uses the fact that $B(\ell,n)$ is the number of subgroups of $\mathbb{Z}^{\ell}$ of index $n$ (a remark by Stanley already mentioned in~\cite{BryanF}) and the main part of the argument is the computation of this number using Hermite normal forms, i.e., Gaussian elimination over the integers. Note that $B(\ell,n)$ is a well studied quantity, see, e.g.,~\cite[Ch. 15]{LubotzkyS} as well as the article by Solomon~\cite{Solomon} where work on $B(\ell,n)$ is traced back to the time of Hermite and Eisenstein.
Also note that a proof of the $x=1$ evaluation of the $\ell=3$ case of (\ref{BFidentity}) was also given in~\cite{Britnell}.
Our proof, given in the next section, is elementary and in the spirit of bijective enumerative combinatorics. In Lemma \ref{polymerlem}, we
 reduce the $A(\ell,n,k)$ to the $k=1$ case of transitive actions, via a polymer gas representation~\cite{Brydges}, in the language of statistical mechanics, or the exponential formula in enumerative combinatorics, often mentioned as the general slogan ``sums over all objects are exponentials of sums over connected objects''. The main argument is a reduction of $A(\ell+1,n,1)$ to the computation of $A(\ell,n,1)$. We condition the sum over tuples $(\sigma_1,\ldots,\sigma_{\ell+1})$, first on the number $r$ of orbits for the sub-tuple $(\sigma_1,\ldots,\sigma_{\ell})$ and then on the set partition $X=\{X_1,\ldots,X_r\}$ of $[n]$ given by that orbit decomposition. With $r$ and $X$ fixed, we then construct a bijection
\begin{equation}
(\sigma_1,\ldots,\sigma_{\ell+1})\longmapsto (\widetilde{\sigma},\gamma,\tau,z)
\label{bijeq}
\end{equation}
where $\widetilde{\sigma}$ is a transitive $\ell$-tuple of commuting permutations on the subset $X_1$ containing the element $1\in[n]$. By $\gamma$ we denote a permutation of $[r]$ which is such that $\gamma(1)=1$. The $\tau$ is a certain collection of bijective maps between blocks $X_i$. Finally, the crucial ingredient is $z$ which is an element of $X_1$.
One can intuitively understand our proof as counting possibly flat or degenerate discrete $(\ell+1)$-dimensional tori with $n$ points. As is familiar in topology, one can build such a torus by gluing both ends of a cylinder. However, we are allowed to perform a twist when doing this gluing and this is determined by $z$. Namely, $(\sigma_{\ell+1})^{r}$, the ``Poincar\'e return map'' to $X_1$, does not necessarily fix $1$ but may send it to some $z\neq 1$.
We remark that it is possible to explicitly iterate the bijection involved in the $\ell+1$ to $\ell$ reduction, but given the complexity of the resulting recursive combinatorial data, we will refrain from doing this here.

\section{Proofs}

We first take care of the reduction to the transitive action case.

\begin{lemma}
We have
\[
A(\ell,n,k)=\ \frac{n!}{k!}\times
\sum_{n_1,\ldots,n_k\ge 1}
\bbone \{n_1+\cdots+n_k=n\}\times \prod_{i=1}^{k}\left(\frac{A(\ell,n_i,1)}{n_i!} \right)\ . 
\]
\label{polymerlem}
\end{lemma}

\noindent{\bf Proof:}
For a tuple $(\sigma_1,\ldots,\sigma_{\ell})$ in $\mathscr{C}_{\ell,n,k}$, let $\Pi(\sigma_1,\ldots,\sigma_{\ell})$ denote the unordered set partition of $[n]$ given by the orbits of the action of the subgroup $\langle\sigma_1,\ldots,\sigma_{\ell}\rangle$. We condition the sum over tuples in $\mathscr{C}_{\ell,n,k}$, according to this set partition. We also sum over orderings of the blocks of that partition (with $k$ blocks), and compensate for this overcounting by dividing by $k!$. This gives
\[
A(\ell,n,k)=\frac{1}{k!}\times\sum_{(X_1,\ldots,X_k)}\ \ \ 
\sum_{(\sigma_1,\ldots,\sigma_{\ell})\in\mathscr{C}_{\ell,n,k}}\ 
\bbone\left\{ \ \Pi(\sigma_1,\ldots,\sigma_{\ell})=\{X_1,\ldots,X_k\}\ \right\}\ ,
\]
where the sum is over ordered tuples of subsets $(X_1,\ldots,X_r)$, where the $X_i$ are nonempty, pairwise disjoint, and together have union equal to $[n]$.
For $1\le i\le k$ and $1\le j\le \ell$, we let $\sigma_j^{(i)}$ be the restriction and corestriction of $\sigma_j$ to the subset $X_i$ which must be stable by $\sigma_j$.
For fixed $X_1,\ldots,X_k$, the sum over tuples $(\sigma_1,\ldots,\sigma_{\ell})$ clearly amounts to summing independently over the tuples $(\sigma_{1}^{(i)},\ldots\sigma_{\ell}^{(i)})$ in each $X_i$, $1\le i\le k$. The tuple $(\sigma_{1}^{(i)},\ldots\sigma_{\ell}^{(i)})$ is made of commuting permutations of $X_i$ whose action on the latter must be transitive. The number of such tuples only depends on the size $|X_i|$ of the set $X_i$, and not its location within $[n]$. As a result, we have
\begin{eqnarray*}
A(\ell,n,k) &=&  \frac{1}{k!}\times
 \sum_{(X_1,\ldots,X_k)} A(\ell,|X_1|,1)\cdots A(\ell,|X_k|,1)\\
 & = & \frac{1}{k!}\times \sum_{n_1,\ldots,n_k\ge 1}
\bbone \{n_1+\cdots+n_k=n\}\times \frac{n!}{n_1!\cdots n_k!}\times \prod_{i=1}^{k}A(\ell,n_i,1) \ ,
\end{eqnarray*}
where the multinomial coefficient accounts for the number of tuples of disjoint sets $(X_1,\ldots,X_k)$ with fixed cardinalities $n_1,\ldots,n_k$.
\qed

We now move on to the main part of the proof, i.e., the $\ell+1$ to $\ell$ reduction and showing that
\begin{equation}
A(\ell+1,n,1)=\sum_{rs=n} A(\ell,s,1)
\times \frac{n!}{r!\times s!^r}\times
(r-1)!\times s!^{r-1}\times s\ ,
\label{reductioneq}
\end{equation}
where the sum is over pairs of integers $r,s\ge 1$ whose product is $n$.
Let $(\sigma_1,\ldots,\sigma_{\ell+1})\in\mathscr{C}_{\ell+1,n,1}$ denote a $(\ell+1)$-tuple of commuting permutations being counted on the left-hand side. We let $X=\{X_1,\ldots,X_r\}:=\Pi(\sigma_1,\ldots,\sigma_{\ell})$ be the set of orbits determined by the first $\ell$ permutations. For a fixed set partition $X$ of $[n]$, define $\mathscr{C}_{\ell+1,n,1}^X\subset\mathscr{C}_{\ell+1,n,1}$ as the set of $(\ell+1)$-tuples which produce the given $X$ by the above definition. We organize the count by conditioning on $X$, i.e., writing
\[
A(\ell+1,n,1)=\sum_X \left|\mathscr{C}_{\ell+1,n,1}^X\right|\ ,
\]
and then computing the terms in the last sum by constructing an explicit bijection between $\mathscr{C}_{\ell+1,n,1}^X$ and a set of combinatorial data whose cardinality is easy to derive.
We will use an automatic numbering of the blocks of $X$ by ordering them according to their minimal element, with respect to the ordered set $[n]$. We let $X_1$ be the block containing the element $1\in[n]$, and number the other blocks so that
\[
1=\min X_1<\min X_2<\cdots<\min X_r\ .
\]

\begin{lemma}
Let $f$ be an element of $\langle\sigma_{\ell+1}\rangle$, i.e., a power of $\sigma_{\ell+1}$, and let $\alpha,\beta\in[r]$.
If $\exists x\in X_{\alpha}$, $f(x)\in X_{\beta}$, then $\forall y\in X_{\alpha}$, $f(y)\in X_{\beta}$.
\end{lemma}

\noindent{\bf Proof:}
Since such $y$ is in the same $\langle\sigma_1,\ldots,\sigma_{\ell}\rangle$-orbit as $x$, there exists a permutation $g\in\langle\sigma_1,\ldots,\sigma_{\ell}\rangle$, such that $y=g(x)$. Since $\sigma_1,\ldots,\sigma_{\ell+1}$ commute, then $g$ must commute with $f$, and therefore $f(y)=f(g(x))=g(f(x))$. This shows that $f(y)$ is in the same 
$\langle\sigma_1,\ldots,\sigma_{\ell}\rangle$-orbit as $f(x)$, namely, $X_{\beta}$.
\qed

The last lemma allows us, from an $f\in\langle\sigma_{\ell+1}\rangle$, to construct a map $\widehat{f}:[r]\rightarrow[r]$
defined by $\widehat{f}(\alpha)=\beta$, whenever $\exists x\in X_{\alpha}$, $f(x)\in X_{\beta}$. 
This construction satisfies $\widehat{{\rm Id}}={\rm Id}$, and $\widehat{f\circ g}=\widehat{f}\circ\widehat{g}$, namely, it gives a group homomorphism from $\langle\sigma_{\ell+1}\rangle$ to $\mathfrak{S}_r$.
We apply this to $f=\sigma_{\ell+1}$ and consider the cycle decomposition of the permutation $\widehat{\sigma_{\ell+1}}$, and focus on the cycle containing the element $1\in[r]$, namely $(\alpha_1\ \alpha_2\ \cdots\ \alpha_t)$, with $\alpha_1=1$.
We clearly have
\[
\sigma_{\ell+1}(X_1)\subset X_{\alpha_2}\ ,\ 
\sigma_{\ell+1}(X_{\alpha_2})\subset X_{\alpha_3}\ ,\ \cdots\ ,\ 
\sigma_{\ell+1}(X_{\alpha_{t-1}})\subset X_{\alpha_t}\ ,\ 
\sigma_{\ell+1}(X_{\alpha_t})\subset X_{1}\ .
\]
Hence $X_{1}\cup X_{\alpha_2}\cup\cdots\cup X_{\alpha_t}$ is stable by $\sigma_{\ell+1}$, in addition to being stable by $\langle\sigma_1,\ldots,\sigma_{\ell}\rangle$ since, each of the $X$ blocks are. Given that the $(\ell+1)$-tuple of
permutations $(\sigma_1,\ldots,\sigma_{\ell+1})$ is assumed to act transitively, this can only happen if the previous union of $X$ blocks is all of $[n]$, i.e., if $t=r$.
For notational convenience, we define the permutation $\gamma\in\mathfrak{S}_r$, by letting $\gamma(i)=\alpha_i$ for all $i\in[r]$. In particular, $\gamma(1)=1$, by construction.
We now have,
\begin{equation}
\sigma_{\ell+1}(X_1)\subset X_{\gamma(2)}\ ,\ 
\sigma_{\ell+1}(X_{\gamma(2)})\subset X_{\gamma(3)}\ ,\ \cdots\ ,\ 
\sigma_{\ell+1}(X_{\gamma(r-1)})\subset X_{\gamma(r)}\ ,\ 
\sigma_{\ell+1}(X_{\gamma(r)})\subset X_{1}\ .
\label{cyclicincleq}
\end{equation}
Since $\sigma_{\ell+1}$ is injective, it follows that
\[
|X_1|\le |X_{\gamma(2)}|\le\cdots\le|X_{\gamma(r)}|\le|X_1|\ ,
\]
and, therefore, all the $X$ blocks must have the same cardinality say $s$, so that $n=rs$, namely, $r$ must divide $n$.
The above argument also produces bijective maps
\[
\tau_{i}:X_{\gamma(i)}\longrightarrow X_{\gamma(i+1)}\ ,
\]
for $1\le i\le r-1$, obtained by restriction (and corestriction) of $\sigma_{\ell+1}$. We collect them into a tuple $\tau=(\tau_1,\ldots,\tau_{r-1})$.
We now define the $\ell$-tuple of permutations of the first block $X_1$ given by
$\widetilde{\sigma}=(\widetilde{\sigma}_1,\ldots,\widetilde{\sigma}_{\ell})$ where, for all $j\in[\ell]$, $\widetilde{\sigma}_j$ is obtained from $\sigma_j$ by restricting it to the subset $X_1$. It is easy to see that $\widetilde{\sigma}$ is an
$\ell$-tuple of commuting permutations of the set $X_1$, which altogether act transitively on it.
Finally, we define the element $z=(\sigma_{\ell+1})^r(1)$
of the block $X_1$. This concludes the definition of the map mentioned in (\ref{bijeq})
which to a tuple $(\sigma_1,\ldots,\sigma_{\ell+1})\in\mathscr{C}_{\ell+1,n,1}$ associates the data $(\widetilde{\sigma},\gamma,\tau,z)$.
Once we establish that this construction is bijective, the reduction formula (\ref{reductioneq})
will follow easily.
Indeed, after identification of $X_1$ with $[s]$, we see that there are $A(\ell,s,1)$ possible choices for $\widetilde{\sigma}$. Deciding on the permutation $\gamma$, which fixes $1$, results in $(r-1)!$ choices. The number of possibilities for the bijective maps in $\tau$ accounts for a factor $s!^{r-1}$, and there are $s$ possibilities for $z$.
Summing over the unordered set partition $X$ can be done with the multinomial coefficient $n!/s!^r$ for ordered set partitions and correcting for the overcounting by dividing by $r!$, as in the proof of Lemma \ref{polymerlem}. All that remains in order to complete the proof of (\ref{reductioneq}) is to show our map (\ref{bijeq}) is indeed bijective. 

\medskip\noindent
{\bf Injectivity:}
We will show how the tuple $(\sigma_1,\ldots,\sigma_{\ell+1})$ is determined by the data $(\widetilde{\sigma},\gamma,\tau,z)$, and the a priori knowledge of the fixed partition $X$. By construction, for all $j$, $1\le j\le \ell$, the restriction of $\sigma_j$ to $X_1$ must be
\begin{equation}
\sigma_j|_{X_1}=\widetilde{\sigma}_j\ .
\label{jX1def}
\end{equation}
Strictly speaking, there is also a change of codomain involved (from $X_1$ to $[n]$), but we ignored this and will continue to do this for the next similar statements.
We must also have, for all $i$, $1\le i\le r-1$,
\begin{equation}
\sigma_{\ell+1}|_{X_{\gamma(i)}}=\tau_i\ .
\label{pplus1Xidef}
\end{equation}
From the commutation relation $\sigma_j\circ(\sigma_{\ell+1})^i=(\sigma_{\ell+1})^i\circ\sigma_j$, restricted to $X_1$,
we deduce that for all $i$, $2\le i\le r$, we must have
\[
\sigma_j\circ\tau_{i-1}\circ\cdots\circ\tau_1=
\tau_{i-1}\circ\cdots\circ\tau_1\circ\widetilde{\sigma}_j
\]
i.e.,
\begin{equation}
\sigma_j|_{X_{\gamma(i)}}=\tau_{i-1}\circ\cdots\circ\tau_1\circ\widetilde{\sigma}_j
\circ\tau_{1}^{-1}\circ\cdots\circ\tau_{i-1}^{-1}\ .
\label{jXidef}
\end{equation}
Hence $\sigma_1,\ldots,\sigma_{\ell}$ are known, while $\sigma_{\ell+1}$
is almost entirely determined. We are only missing the restriction of $\sigma_{\ell+1}$ on the last block $X_{\gamma(r)}$.
Since $z$ is in the orbit $X_1$ of the element $1$ for the action of $\sigma_1,\ldots,\sigma_{\ell}$, or equivalently $\widetilde{\sigma}_1,\ldots,\widetilde{\sigma}_{\ell}$, 
there exists $g\in\langle\widetilde{\sigma}_1,\ldots,\widetilde{\sigma}_{\ell}\rangle$, such that $g(1)=z$. We claim that we must have
\begin{equation}
\sigma_{\ell+1}|_{X_{\gamma(r)}}=g\circ\tau_{1}^{-1}\circ\cdots\circ\tau_{r-1}^{-1}\ .
\label{pplus1Xrdef}
\end{equation}
Indeed, let $x\in X_{\gamma(r)}$, then $x=(\sigma_{\ell+1})^{r-1}(y)$ for some $y\in X_1$. Again, by transitivity on $X_1$, there exists $h\in\langle\sigma_1,\ldots,\sigma_{\ell}\rangle$ such that $y=h(1)$. As a consequence of the Abelian property of the group $\langle\sigma_1,\ldots,\sigma_{\ell+1}\rangle$, we must have
\begin{eqnarray*}
\sigma_{\ell+1}(x) &=& (\sigma_{\ell+1})^r\circ h(1)\\
 &=& h\circ(\sigma_{\ell+1})^r(1)\\
&= & h(z)\\
& = & h(g(1))\\
&=& g(h(1))\\
& = & g(y)\\
& = & g\circ\tau_{1}^{-1}\circ\cdots\circ\tau_{r-1}^{-1}(x)\ .
\end{eqnarray*}
We now have recovered the restrictions of all $\ell+1$ permutations $\sigma_j$ on all 
blocks $X_i$ of the decomposition of $[n]$, from the output of our map, which shows that it is injective.

\medskip\noindent
{\bf Surjectivity:}
We start from the data $(\widetilde{\sigma},\gamma,\tau,z)$
and construct $(\sigma_1,\ldots,\sigma_{\ell+1})\in\mathscr{C}_{\ell+1,n,1}^{X}$ which maps to it.
This time, we use the equations (\ref{jX1def}), (\ref{pplus1Xidef}), (\ref{jXidef}), (\ref{pplus1Xrdef}) as definitions of $\sigma_1,\ldots,\sigma_{\ell+1}$ as maps $[n]\rightarrow[n]$. The use of (\ref{pplus1Xrdef}) requires some care, namely showing the uniqueness of $g$.
Let $\widetilde{H}:=\langle\widetilde{\sigma}_1,\ldots,\widetilde{\sigma}_p\rangle$. The hypothesis on the tuple $\widetilde{\sigma}$ is that it is made of $\ell$ commuting permutations of the set $X_1$, such that the permutation action of $\widetilde{H}$ on $X_1$ is transitive. Suppose $g_1(1)=g_2(1)=z$ for some $g_1,g_2\in\widetilde{H}$.
If $x\in X_1$, then $\exists h\in\widetilde{H}$, $h(1)=x$. By the Abelian property of $\widetilde{H}$, we have
\[
g_i(x)=g_i\circ h(1)=h\circ g_i(1)=h(z)\ ,
\]
for $i=1$ as well as $i=2$, and thus $g_1(x)=g_2(x)$. Since $x$ is arbitrary, we have $g_1=g_2$. This justifies the use of (\ref{pplus1Xrdef}) as a definition of a map.
We now have constructed the maps $\sigma_1,\ldots,\sigma_{\ell+1}$.
It is immediate, from (\ref{jX1def}) and (\ref{jXidef}), that $\sigma_1,\ldots,\sigma_{\ell}$ are bijective within each $X_{\gamma(i)}$, $1\le i\le r$, and therefore over all of $[n]$.
One easily checks also the commutation relations $\sigma_j\circ\sigma_{j'}=\sigma_{j'}\circ\sigma_j$, $1\le j,j'\le \ell$, on each $X$ block, and therefore on $[n]$.
From (\ref{pplus1Xidef}), we see that $\sigma_{\ell+1}$ is injective on each $X_{\gamma(i)}$, $1\le i\le r-1$, and the images of these restrictions are disjoint because $\gamma$ is a permutation. From (\ref{pplus1Xrdef}),
it holds that $\sigma_{\ell+1}|_{X_{\gamma(r)}}:X_{\gamma(r)}\rightarrow X_1$ is bijective. As a result, $\sigma_{\ell+1}:[n]\rightarrow[n]$ is bijective. 
From (\ref{pplus1Xidef}) and (\ref{jXidef}), we also obtain
\[
\sigma_j\circ\sigma_{\ell+1}|_{X_{\gamma(i)}}=
\tau_{i}\circ\cdots\circ\tau_1\circ\widetilde{\sigma}_j
\circ\tau_{1}^{-1}\circ\cdots\circ\tau_{i-1}^{-1}
=\sigma_{\ell+1}\circ\sigma_{j}|_{X_{\gamma(i)}}\ ,
\]
for all $i,j$ such that $1\le j\le \ell$ and $1\le i\le r-1$.
Finally, for all $j$, $1\le j\le \ell$, the restrictions of $\sigma_j\circ\sigma_{\ell+1}$ and $\sigma_{\ell+1}\circ\sigma_{j}$ on $X_{\gamma(r)}$ coincide, because $g$ and $\widetilde{\sigma}_j$ must commute. We have now checked that $(\sigma_1,\ldots,\sigma_{\ell+1})$ is a commuting tuple of permutations of $[n]$.
The corresponding action is transitive because (\ref{cyclicincleq}) holds by construction and $\widetilde{\sigma}$ is assumed to act transitively on $X_1$.
Checking that the produced tuple
$(\sigma_1,\ldots,\sigma_{\ell+1})\in\mathscr{C}_{\ell+1,n,1}^{X}$ indeed maps to
$(\widetilde{\sigma},\gamma,\tau,z)$ is straightforward. Therefore, our map is surjective.

\medskip
In order to finish the proof of Theorem \ref{mainthm}, we define $C(\ell,n):=\frac{A(\ell,n,1)}{(n-1)!}$. Since $A(1,n,1)=(n-1)!$ counts cyclic permutations of $n$ elements, we have $C(1,n)=1=B(1,n)$. The, now established, recursion (\ref{reductioneq}) implies that $C$ satisfies
\[
C(\ell+1,n)=\sum_{rs=n}s\ C(\ell,s)\ .
\]
By a trivial induction on $\ell$, $C(\ell,n)$ must coincide with $B(\ell,n)$ defined, e.g., in (\ref{altdefeq}). We plug $A(\ell,n,1)=(n-1)!\times B(\ell,n)$ in the result of Lemma \ref{polymerlem}, and Theorem \ref{mainthm} follows. 
\qed

\section{On conjecture \ref{mainconj}}\label{conjsec}
\label{conjdiscsec}

As mentioned in the introduction, the case $\ell=1$ of Conjecture \ref{mainconj} is well established. The opposite extreme ``$\ell=\infty$'' is settled in the companion article~\cite{Abdesselam2}.
Let us now focus on the $\ell=2$ case, and relate it to an already large body of literature, in particular, the work of Heim, Neuhauser, and many others.
Since, for $\ell=2$, $B(\ell-1,m)=B(1,m)=1$, the Bryan-Fulman identity (\ref{BFidentity})
simply reads
\[
\sum_{n=0}^{\infty}\sum_{k=0}^{n}\frac{1}{n!}A(2,n,k)x^k u^n=
\prod_{m=1}^{\infty}
(1-u^{m})^{-x}\ .
\]
On the other hand, the so-called D'Arcais polynomials $P_n(x)$ are defined~\cite{DArcais} by the generating function identity
\[
\prod_{m=1}^{\infty}
(1-u^{m})^{-x}=\sum_{n=0}^{\infty} P_n(x) u^n\ . 
\]
The D'Arcais polynomials can therefore be expressed in terms of commuting pairs of permutations
\begin{equation}
P_n(x)=\frac{1}{n!}\sum_{k=0}^n A(2,n,k)\ x^k\ .
\label{DArcaiseq}
\end{equation}
We are not aware of the commuting permutation interpretation (\ref{DArcaiseq}) of D'Arcais polynomials having been used in the number theory literature reviewed, e.g., in~\cite{HeimN}, and we hope it could be of help in this area.
If one shifts the variable $x$ by one, one gets the standard formulation of the Nekrasov-Okounkov formula~\cite{NekrasovO,Westbury}
\[
\prod_{m=1}^{\infty}
(1-u^{m})^{-x-1}=\sum_{n=0}^{\infty} Q_n(x) u^n
\]
where
\[
Q_n(x)=\sum_{\lambda\vdash n}\prod_{\Box\in\lambda}
\left(1+\frac{x}{h(\Box)^2}\right)\ .
\]
Namely, the sum is over integer partitions $\lambda$ of $n$. The product is over cells in the usual Ferrers-Young diagram of the partition $\lambda$, and $h(\Box)$ denotes the hook length number of that cell.
Clearly $Q_n(x)=P_n(x+1)$ and therefore, the log-concavity (of the coefficients of) the polynomial $P_n$ would imply that of $Q_n$ as well as the unimodality of the latter which was conjectured by Heim and Neuhauser as well as Amdeberhan (see~\cite{HeimN} and references therein). As a strengthening of this unimodality conjecture, the log-concavity of the $P_n(x)$'s, i.e., the $\ell=2$ case of Conjecture \ref{mainconj} was stated as Challenge 3 in~\cite{HeimN}. The authors also reported on checking this numerically for all $n\le 1500$.
While the log-concavity in the $\ell=1$ case can be derived using the real-rootedness of the relevant polynomial, namely, the Pochhammer symbol, this approach cannot work for $\ell=2$. Indeed, D'Arcais polynomials can have roots off the real axis, as was shown in~\cite{HeimN2}.
For recent progress towards such log-concavity properties in the $\ell=2$ case, see~\cite{HongZ,Zhang}.

\medskip
Using Mathematica, we checked that Conjecture \ref{mainconj} is true for $\ell=3,4,5$ for all $n\le 100$. One can also test the conjecture by considering the dilute polymer gas regime, in the terminology of statistical mechanics (see, e.g.,~\cite{Brydges}), i.e., when $k$ is close to $n$ and most orbits are singletons, as in the next proposition. Note that the latter can also be deduced from~\cite[Proposition 4]{HeimNarXiv}. The $\ell=2$ case was, in fact, already proved in~\cite[Corollary 4]{HeimNarXiv}.

\begin{proposition} 
The inequality in Conjecture \ref{mainconj} holds for all $\ell\ge 1$, and $n\ge 3$, when $k=n-1$.
\end{proposition}

\noindent{\bf Proof:}
Let
\[
\Delta(\ell,n):=A(\ell,n,n-1)^2-A(\ell,n,n)\ A(\ell,n,n-2)\ .
\]
From Theorem \ref{mainthm}, we easily deduce
\begin{eqnarray*}
A(\ell,n,n) & = & 1 \\
A(\ell,n,n-1) & = & \binom{n}{2}\ (2^{\ell}-1) \\
A(\ell,n,n-2) & = & \binom{n}{3}\ (3^{\ell}-1)+\binom{n}{4}\ 3 (2^{\ell}-1)^2\ .
\end{eqnarray*}
Therefore
\[
\Delta(\ell,n)=\left[{\binom{n}{2}}^2-3\binom{n}{4}\right] (2^{\ell}-1)^2-\binom{n}{3}\ (3^{\ell}-1)\ .
\]
As mentioned before, the conjecture is known for $\ell=1$, so now we focus on $\ell\ge 2$.
If $\ell\ge 3$, then $2\left(\frac{1}{2}\right)^{\ell}+\left(\frac{3}{4}\right)^{\ell}\le \frac{43}{64}$, the $\ell=3$ value. Therefore, for $\ell\ge 3$, we have $4^{\ell}\ge 2\times 2^{\ell}+3^{\ell}$ which implies
\[
4^{\ell}-2\times 2^{\ell}+1\ge 3^{\ell}-1\ .
\]
The last inequality being also true for $\ell=2$, we have that for all $\ell\ge 2$, the inequality $(2^{\ell}-1)^2\ge 3^{\ell}-1$ holds.
Hence
\begin{eqnarray*}
\Delta(\ell,n) &\ge & \left[{\binom{n}{2}}^2-3\binom{n}{4}-\binom{n}{3}\right]
(2^{\ell}-1)^2 \\
 & = & \frac{1}{24}n(n-1)(3n^2+5n-10)(2^{\ell}-1)^2\ .
\end{eqnarray*}
Since $n\ge 3$ implies $3n^2+5n-10\ge 32>0$, we have $\Delta(\ell,n)>0$.
\qed

\bigskip
\noindent{\bf Acknowledgements:}
{\small
The first author thanks Ken Ono for introducing him to the Nekrasov-Okounkov formula, and the unimodality conjecture of Amdeberhan, Heim and Neuhauser. We also thank Bernhard Heim, Markus Neuhauser, and Alexander Thomas for a careful reading of the first version of this article. We also thank the anonymous referee for useful comments that helped improve this article.
}


\bigskip
\noindent{\bf Conflict of interest statement:} 
{\small
On behalf of all the authors, the corresponding author states that there is no conflict of interest. 
}


\bigskip
\noindent{\bf Data availability statement:}
{\small
The Mathematica files for the computations presented in \S\ref{conjdiscsec} are available from the corresponding author, upon request.
}


\end{document}